\documentclass[12pt,leqno]{article}
\usepackage{amsfonts}
\usepackage{amsmath, amsthm, amsfonts, amssymb, color}
\usepackage{mathrsfs}
\usepackage{url}
\usepackage{color}
\theoremstyle{definition}

\newcommand{\scr}[1]{\mathscr #1}
\definecolor{wco}{rgb}{0.5,0.2,0.3}

\numberwithin{equation}{section} \theoremstyle{remark}

\newcommand{\ua}{\uparrow}

\title{{\bf Precise Limit in  Wasserstein Distance for Conditional Empirical Measures of Dirichlet Diffusion Processes}\footnote{Supported in
 part by  the National Key R\&D Program of China (No. 2020YFA0712900),  NNSFC (11771326, 11831014, 11921001), and DFG through the CRC Taming uncertainty and profiting from randomness
and low regularity in analysis, stochastics and their applications.} }
\author{
{\bf    Feng-Yu Wang$^{a),b)}$    }\\
\footnotesize{$^{a)}$ Center for Applied Mathematics, Tianjin University, Tianjin 300072, China }\\
 \footnotesize{ $^{b)}$ Department of Mathematics,
Swansea University,
Bay Campus,
Swansea,
SA1 8EN, United Kingdom}  }

\begin{document}
\allowdisplaybreaks
\def\R{\mathbb R}  \def\ff{\frac} \def\ss{\sqrt} \def\B{\mathbf
B}\def\TO{\mathbb T}
\def\I{\mathbb I_{\pp M}}\def\p<{\preceq}
\def\N{\mathbb N} \def\kk{\kappa} \def\m{{\bf m}}
\def\ee{\varepsilon}\def\ddd{D^*}
\def\dd{\delta} \def\DD{\Delta} \def\vv{\varepsilon} \def\rr{\rho}
\def\<{\langle} \def\>{\rangle} \def\GG{\Gamma} \def\gg{\gamma}
  \def\nn{\nabla} \def\pp{\partial} \def\E{\mathbb E}
\def\d{\text{\rm{d}}} \def\bb{\beta} \def\aa{\alpha} \def\D{\scr D}
  \def\si{\sigma} \def\ess{\text{\rm{ess}}}
\def\beg{\begin} \def\beq{\begin{equation}}  \def\F{\scr F}
\def\Ric{{\rm Ric}} \def\Hess{\text{\rm{Hess}}}
\def\e{\text{\rm{e}}} \def\ua{\underline a} \def\OO{\Omega}  \def\oo{\omega}
 \def\tt{\tilde}
\def\cut{\text{\rm{cut}}} \def\P{\mathbb P} \def\ifn{I_n(f^{\bigotimes n})}
\def\C{\scr C}      \def\aaa{\mathbf{r}}     \def\r{r}
\def\gap{\text{\rm{gap}}} \def\prr{\pi_{{\bf m},\varrho}}  \def\r{\mathbf r}
\def\Z{\mathbb Z} \def\vrr{\varrho} \def\ll{\lambda}
\def\L{\scr L}\def\Tt{\tt} \def\TT{\tt}\def\II{\mathbb I}
\def\i{{\rm in}}\def\Sect{{\rm Sect}}  \def\H{\mathbb H}
\def\M{\scr M}\def\Q{\mathbb Q} \def\texto{\text{o}} \def\LL{\Lambda}
\def\Rank{{\rm Rank}} \def\B{\scr B} \def\i{{\rm i}} \def\HR{\hat{\R}^d}
\def\to{\rightarrow}\def\l{\ell}\def\iint{\int}
\def\EE{\scr E}\def\Cut{{\rm Cut}}\def\W{\mathbb W}
\def\A{\scr A} \def\Lip{{\rm Lip}}\def\S{\mathbb S}
\def\BB{\scr B}\def\Ent{{\rm Ent}} \def\i{{\rm i}}\def\itparallel{{\it\parallel}}
\def\g{{\mathbf g}}\def\Sect{{\mathcal Sec}}\def\T{\mathcal T}\def\V{{\bf V}}
\def\PP{{\bf P}}\def\HL{{\bf L}}\def\Id{{\rm Id}}\def\f{{\bf f}}\def\cut{{\rm cut}}
\def\sm{\preceq}

\def\BL{\scr A}

\maketitle

\begin{abstract} Let $M$ be a   $d$-dimensional connected compact Riemannian manifold   with      boundary $\partial M$, let $V\in C^2(M)$ such that $\mu({\rm d}x):={\rm e}^{V(x)}{\rm d} x$ is a probability measure,  and let
$X_t$ be the diffusion process generated by $L:=\Delta+\nabla V$ with  $\tau:=\inf\{t\ge 0: X_t\in\partial M\}$. Consider the conditional empirical measure
 $\mu_t^\nu:= \mathbb E^\nu\big(\frac 1 t \int_0^t \delta_{X_s}{\rm d} s\big|t<\tau\big)$ for the diffusion process    with initial  distribution $\nu$ such that  $\nu(\partial M)<1$. Then
$$\lim_{t\to\infty} \big\{t\mathbb W_2(\mu_t^\nu,\mu_0)\big\}^2 =  \frac 1 {\{\mu(\phi_0)\nu(\phi_0)\}^2} \sum_{m=1}^\infty \frac{\{\nu(\phi_0)\mu(\phi_m)+ \mu(\phi_0) \nu(\phi_m)\}^2}{(\lambda_m-\lambda_0)^3},$$
where $\nu(f):=\int_Mf {\rm d} \nu$ for a measure $\nu$ and $f\in L^1(\nu)$, $\mu_0:=\phi_0^2\mu$, $\{\phi_m\}_{m\ge 0}$ is the eigenbasis of $-L$   in $L^2(\mu)$ with the Dirichlet boundary,
  $\{\lambda_m\}_{m\ge 0}$ are the corresponding Dirichlet eigenvalues,
  and $\mathbb W_2$ is the $L^2$-Wasserstein distance induced by the Riemannian metric.
 \end{abstract} \noindent
 \vskip 2cm

\section{Introduction }

Let $M$ be a $d$-dimensional    connected  compact Riemannian manifold    with a smooth  boundary $\pp M$.  Let $V\in C^2(M)$ such that $\mu(\d x)=\e^{V(x)}\d x$ is a probability measure on $M$, where $\d x$ is the Riemannian volume measure. 
Let $X_t$  be the diffusion process generated by $L:=\DD+\nn V$ with hitting time
$$\tau:=\inf\{t\ge 0: X_t\in\pp M\}.$$ Here, according to the convention in Riemannian geometry,  the vector field $\nn V$ is regarded as a first-order differential operator with
$(\nn V)f:= \<\nn V,\nn f\>$ for differentiable functions $f$.
Denote by $\scr P$   the set of all probability measures on $M$, and let $\E^\nu$ be the expectation taken for the diffusion process with initial distribution $\nu\in\scr P$. 
 Consider the conditional empirical measure
$$\mu_t^\nu:=\E^\nu\bigg(\ff 1 t \int_0^t \dd_{X_s}\d s\bigg|t<\tau\bigg),\ \ t>0,\nu\in\scr P.$$
Since $\tau=0$ when $X_0\in \pp M$, to ensure $\P^\nu(\tau>t)>0$ we only consider
$$\nu\in \scr P_0:=\big\{\nu\in \scr P:\ \nu(M^\circ)>0\big\},\ \ M^\circ:= M\setminus \pp M.$$

Let $\{\phi_m\}_{m\ge 0}$ be  the eigenbasis in $L^2(\mu)$ of $-L$   with the Dirichlet boundary such that $\phi_0>0$ in $M^\circ$,   and let $\{\ll_m\}_{m\ge 0}$ be the associated eigenvalues listed in the increasing order counting multiplicities; that is, $\{\phi_m\}_{m\ge 0}$ is an orthonormal basis of    $L^2(\mu)$  such that
$$L \phi_m= -\ll_m \phi_m,\ \ m\ge 0.$$
Then $\mu_0:=\phi_0^2\mu$ is a probability measure on $M$.
It is easy to see from \cite[Theorem 2.1]{CJ14} that for any probability measure $\nu$ supported on $M^\circ$, we have
$$\lim_{t\to\infty} \|\mu_t^\nu-\mu_0\|_{var}=0,$$
where $\|\cdot\|_{var}$ is the total variational norm.

In this paper, we investigate the convergence of $\mu_t^\nu$ to $\mu_0$
under the Wasserstein distance $\W_2$:
 $$\W_2(\mu_1,\mu_2):= \inf_{\pi\in \C(\mu_1,\mu_2)} \bigg(\int_{M\times M} \rr(x,y)^2 \pi(\d x,\d y) \bigg)^{\ff 1 2},\ \ \mu_1,\mu_2\in \scr P,$$
 where $\C(\mu_1,\mu_2)$ is the set of all probability measures on $M\times M$ with marginal distributions $\mu_1$ and $\mu_2$, and  $\rr(x,y)$ is the Riemannian distance between $x$ and $y$, i.e.   the length of the shortest curve on $M$ linking $x$ and $y$.

Recently, the convergence rate under $\W_2$ has been characterized  in \cite{WZ20} for the empirical measures of the $L$-diffusion processes without boundary (i.e. $\pp M=\emptyset$) or with a reflecting boundary.
Since in the present setting  the diffusion process is killed at time $\tau$, it is reasonable to consider the conditional empirical measure $\mu_t^\nu$ given $t<\tau$.
This is a counterpart to the quasi-ergodicity for the convergence of the conditional distribution $\tt \mu_t$ of $X_t$ given $t<\tau$. Unlike in the case without boundary or with a reflecting boundary where both
the distribution and the empirical measure of $X_t$
 converge to the unique invariant probability measure, in the present case the conditional distribution $\tt\mu_t$ of $X_t$ given $t<\tau$ converges to $\tt\mu_0:= \ff{\phi_0}{\mu (\phi_0)}\mu $
 rather than   $\mu_0:=\phi_0^2\mu$,
 and this convergence is called the quasi-ergodicity in the literature, see for instance \cite{CMM}  and references within.

  Let $\nu(f):=\int_Mf\d\nu$ for $\nu\in \scr P$ and $f\in L^1(\nu)$. The main result of this paper is the following.

 \beg{thm}\label{T1.1} For any $\nu\in \scr P_0$,
\beg{align*} \lim_{t\to\infty} \big\{ t^2\W_2(\mu_{t}^\nu, \mu_0)^2\big\}  
=I:=  \ff1 {\{\mu(\phi_0)\nu(\phi_0)\}^2} \sum_{m=1}^\infty \ff{\{\nu(\phi_0)\mu(\phi_m)+ \mu(\phi_0) \nu(\phi_m)\}^2}{(\ll_m-\ll_0)^3}>0. \end{align*} 
 If either $d\le 6$ or  $d\ge 7$ but $\nu=h\mu$ with  $ h\in L^{\ff{2d}{d+6}}(\mu)$,   
 then   $I<\infty$. 
 \end{thm}

\paragraph{Remark 1.1.}  
(1) Let $X_t$ be the (reflecting) diffusion process generated by $L$ on $M$ where $\pp M$ may be empty. We consider the mean empirical measure
$\hat \mu_t^\nu := \E (\ff 1 t \int_0^t \dd_{X_s}\d s),$ where $\nu$ is the initial distribution of $X_t$. Then  
\beq\label{70} \lim_{t\to\infty} \big\{t^2\W_2(\hat \mu_{t}^\nu, \mu_0)^2\big\}=    \sum_{m=1}^\infty \ff{ \{\nu(\phi_m)\}^2}{\ll_m^3}<\infty,\end{equation}
where $\{\phi_m\}_{m\ge 1} $ is the eigenbasis of $-L$ in $L^2(\mu)$ with the Neumann boundary condition if $\pp M$ exists,  $\{\ll_m\}_{m\ge 1}$ are the corresponding non-trivial (Neumann) eigenvalues, and the limit is zero if and only if $\nu=\mu$.
This can be confirmed by   the proof  of Theorem \ref{T1.1}     with $\phi_0=1,\ll_0=0$ and $\mu(\phi_m)=0$ for $m\ge 1.$  In this case, $\mu$ is the unique invariant probability measure of $X_t$,  so that
$\hat\mu_t^\mu=\mu$ for $t\ge 0$ and hence the limit in \eqref{70} is zero for $\nu=\mu.$  However, in the Dirichlet diffusion case, the conditional distribution of $(X_s)_{0\le s\le t}$ given $t<\tau$ is no longer stationary, so that even starting from the limit distribution $\mu_0$ we {\bf do not have }
$\mu_t^{\mu_0}=\mu_0$ for $t>0$. This leads to a non-zero limit in Theorem \ref{T1.1} even for $\nu=\mu_0.$

(2) It is also interesting to investigate the convergence of $\E^\nu(\W_2(\mu_t,\mu_0)^2 |t<\tau)$ for $\mu_t:= \ff 1 t\int_0^t \dd_{X_s}\d s,$ which is the counterpart to the study of \cite{WZ20} where the case without boundary or with a reflecting boundary
is considered. According to \cite{WZ20}, the convergence rate  of $\E^\nu(\W_2(\mu_t,\mu_0)^2 |t<\tau)$ will be at most $t^{-1}$, which is slower than the rate $t^{-2} $ for $\W_2(\mu_t^\nu,\mu_0)^2$ as shown in Theorem \ref{T1.1}, see  \cite{W20b} for details, see also \cite{W20c, W20d} for extensions to diffusion processes on non-compact manifolds and SPDEs.

(3) Let $\nu=h\mu$. It is easy to see that $I<\infty$ is equivalent to $h\in \D((-L)^{-\ff 3 2})$.  By the Sobolev inequality, for any $p\in [1, \ff d 3)$,   there exists a constant $K>0$ such that 
\beq\label{SB} \|(-L)^{-\ff 3 2} f\|_{L^{\ff{dp}{d-3p}}(\mu)}\le K \|f\|_{L^p(\mu)},\ \ f\in L^p(\mu).\end{equation} Taking $p= \ff{2d}{d+6}$ which is large than $1$ when $d\ge 7$,   we see that $h\in L^p(\mu)$ implies $h\in \D((-L)^{-\ff 3 2})$ and hence $I<\infty$.   So,   the sharpness of the Sobolev inequality implies that of the condition $h\in L^{\ff{2d}{d+6}} (\mu).$ 

\

In Section 2, we first recall some well known facts on the Dirichlet semigroup, then  present an upper bound estimate  on $\|\nn (\phi_m\phi_0^{-1})\|_\infty$.
The latter   is non-trivial when $\pp M$ is non-convex, and  should be interesting by itself. With these preparations,  we prove upper and lower bound estimates     in  Sections 3 and 4 respectively. 

 \section{Some preparations}

We first recall some well known facts on the Dirichlet semigroup, see for instances \cite{Chavel, Davies, OUB, W14}.
 Let $\{\phi_m\}_{m\ge 0}$ be the eigenbasis of the Dirichlet operator $L$ in $L^2(\mu)$, with Dirichlet eigenvalues $\{\ll_m\}_{m\ge 0}$ of $-L$ listed in the increasing order counting multiplicities; that is, $\{\phi_m\}_{m\ge 0}$ is an orthonormal basis of    $L^2(\mu)$  such that
$$L \phi_m= -\ll_m \phi_m,\ \ m\ge 0.$$
  For simplicity, we denote $a\sm b$ for two positive functions $a$ and $b$ if $a\le cb$ holds for some constant $c>0$. 
 Then $\ll_0>0$ and
 \beq\label{EG}\|\phi_m\|_\infty  \sm \ss m,\ \   m^{\ff 2 d}\sm \ll_m-\ll_0  \sm m^{\ff 2 d},\ \ m\ge 1.\end{equation}
  Let $\rr_\pp$ be the Riemannian distance function to the boundary $\pp M$. Then $\phi_0^{-1}\rr_\pp$ is bounded such that
 \beq\label{JSS} \|\phi_0^{-1}\|_{L^p(\mu_0)}<\infty,\ \ p\in [1,3).\end{equation}
 The Dirichlet heat kernel has the representation
 $$p_t^D(x,y)=\sum_{m=0}^\infty \e^{-\ll_m t} \phi_m(x)\phi_m(y),\ \ t>0, x,y\in M.$$
 Let $\E^x$ denote the expectation for the $L$-diffusion process starting at point $x$. Then Dirichlet diffusion semigroup generated by $L$ is given by
\beq\label{2.1} \beg{split} &P_t^D f(x):= \E^x[f(X_t)1_{\{t<\tau\}}]=\int_Mp_t^D(x,y)f(y)\mu(\d y)\\
&=\sum_{m=0}^\infty \e^{-\ll_m t} \mu(\phi_mf)\phi_m(x),\ \ t>0, f\in L^2(\mu).\end{split}\end{equation}
We have  
 \beq\label{AC0} \|P_t^D\|_{L^p(\mu)\to L^q(\mu)}:= \sup_{\mu(|f|^p)\le 1} \|P_t^Df\|_{L^q(\mu)}\sm \e^{-\ll_0t} (1\land t)^{-\ff {d (q-p)}{2pq}}, \ \ t>0, q\ge p\ge 1.\end{equation}

Next,  let $L_0= L+2\nn\log \phi_0$. Then $L_0$ is a self-adjoint operator in $L^2(\mu_0)$ with semigroup $P_t^0:=\e^{t L_0}$ satisfying
 \beq\label{PR0} P_t^0f=\e^{\ll_0 t}\phi_0^{-1} P_t^D(f\phi_0),\ \ f\in L^2(\mu_0),\ \ t\ge 0.\end{equation}
So,    $\{\phi_0^{-1}\phi_m\}_{m\ge 0}$ is an eigenbasis of $L_0$ in $L^2(\mu_0)$ with
 \beq\label{PR1} L_0(\phi_m\phi_0^{-1})= -(\ll_m-\ll_0)\phi_m\phi_0^{-1},\ \ P_t^0(\phi_m\phi_0^{-1})=\e^{-(\ll_m-\ll_0)t} \phi_m\phi_0^{-1},\ \ m\ge 0, t\ge 0.\end{equation}
 Consequently,
 \beq\label{ONN} P_t^0f = \sum_{m=0}^\infty \mu_0(f\phi_m\phi_0^{-1}) \e^{-(\ll_m-\ll_0)t} \phi_m\phi_0^{-1},\ \ f\in L^2(\mu_0),\end{equation}
 and the heat kernel of $P_t^0$ with respect to $\mu_0$ is given by
 \beq\label{ON*} p_t^0(x,y)= \sum_{m=0}^\infty (\phi_m\phi_0^{-1})(x)  (\phi_m\phi_0^{-1})(y) \e^{-(\ll_m-\ll_0)t},\ \ x,y\in M, t>0.\end{equation}
By the intrinsic ultracontractivity, see for instance \cite{07OW}, we have 
\beq\label{PRR0}  \|P_t^0-\mu_0\|_{L^1(\mu_0)\to L^\infty(\mu_0)}:= \sup_{\mu_0(|f|)\le 1} \|P_t^0f-\mu_0(f)\|_{\infty} \sm \ff {\e^{-(\ll_1-\ll_0)t}}{(1\land t)^{\ff{d+2}2}},\ \ t>0.\end{equation}
Combining this with the semigroup property and the contraction of $P_t^0$ in $L^p(\mu)$ for any $p\ge 1$, we obtain 
\beq\label{000} \|P_t^0-\mu_0\|_{L^p(\mu_0)}:=  \sup_{\mu_0(|f|^p)\le 1} \|P_t^0f-\mu_0(f)\|_{L^p(\mu_0)} \sm \e^{-(\ll_1-\ll_0)t},\ \ t\ge 0, p\ge 1.\end{equation}
By the interpolation theorem,  \eqref{PRR0} and \eqref{000} yield
\beq\label{**1} \|P_t^0-\mu_0\|_{L^p(\mu_0)\to L^q(\mu_0)}\sm  \e^{-(\ll_1-\ll_0)t} \{1\land t\}^{-\ff{(d+2)(q-p)}{2pq}},\ \ t>0,\infty\ge q>p\ge 1.\end{equation}
Since  $\mu_0(\phi_m^2\phi_0^{-2})=1$,    \eqref{**1} for $p=2$ implies
$$\|\phi_m\phi_0^{-1}\|_\infty = \e^{(\ll_m-\ll_0)t} \|P_t^0(\phi_m\phi_0^{-1})\|_\infty
\sm  \ff{ \e^{(\ll_m-\ll_0)t}}{(1\land t)^{\ff{d+2}4}},\ \ t>0.$$  Taking $t= (\ll_m-\ll_0)^{-1}$ and applying \eqref{EG},
we derive 
\beq\label{CC} \|\phi_m\phi_0^{-1}\|_\infty  \sm m^{\ff{d+2}{2d}},\ \ m\ge 1.\end{equation}

\

In the remainder of this section, we investigate gradient estimates on $P_t^0$ and $\phi_m\phi_0^{-1}$, which will be used in Section 4 for the study of the lower bound estimate on $\W_2(\mu_t^\nu,\mu_0)$.
To this end, we need to estimate the Hessian tensor of $\log \phi_0$.

Let $N$ be the inward unit normal vector field of $\pp M$. We call $M$ (or $\pp M$) convex if
\beq\label{SC} \<\nn_u N, u\>= \Hess_{\rr_\pp}(u,u)\le 0,\ \ u\in T\pp M,\end{equation}
where $\rr_\pp$ is the distance function to the boundary $\pp M$, and $T\pp M$ is the tangent bundle of the $(d-2)$-dimensional manifold $\pp M$.  When $d=1$, the boundary $\pp M$ degenerates to a set of two end points,
such  that $\pp M=\emptyset$ and  the condition \eqref{SC} trivially holds; that is,  $M$ is convex for $d=1$.
Recall that $M^\circ:=M\setminus \pp M$ is the interior of $M$.

\beg{lem}\label{L3.0}   If $\pp M$ is convex, then there exists a constant $K_0\ge 0$ such that
$$\Hess_{\log\phi_0}(u,u)\le K_0|u|^2,\ \ u\in TM^\circ.$$
\end{lem}
\beg{proof} Since $M$ is compact with smooth boundary, there exists a constant $r_0>0$ such that $\rr_\pp$ is smooth on the set
$$\pp_0M:= \{x\in M: \rr_\pp(x)\le r_0\}.$$ Since $\phi_0$ is smooth  and satisfies  $\phi_0\ge c\rr_\pp$ for some constant $c>0$,   we have
$\log(\phi_0\rr_\pp^{-1}) \in C_b^2(\pp_0M)$. So, it suffices to find  a constant $c>0$ such that
\beq\label{HESS} \Hess_{\log\rr_\pp}(u,u)\le c|u|^2,\ \ u\in TM^\circ.\end{equation}
To this end, we fisrt  estimate
$ \Hess_{\rr_\pp}  $ on the boundary $\pp M$.
For any $x\in\pp M$ and $u\in T_x M,$ consider the orthogonal decomposition $u= u_1+u_2$, where
$$u_1= \<N,u\>N, \ \ u_2:=u-u_1\in T\pp M.$$
Since $|\nn\rr_\pp|=1$ on $\pp_0M$, we have
\beq\label{H2} \Hess_{\rr_\pp}(X,N)=\Hess_{\rr_\pp}(X,\nn\rr_\pp)=\ff 1 2  \<X,\nn |\nn\rr_\pp|^2\>=0,\ \ X\in   T_xM.\end{equation}
On the other hand, since $u_2\in T\pp M$ and $\nn\rr_\pp=N$ on $\pp M$,  \eqref{SC} implies
$$\Hess_{\rr_\pp}(u_2,u_2)= \<\nn_{u_2}N,u_2\>\le 0.$$
Combining this with \eqref{H2} we obtain
$$\Hess_{\rr_\pp}(u,u)= \<N, u\>^2\Hess_{\rr_\pp}(N,N)+ 2\<N,u\> \Hess_{\rr_\pp}(u_2, N)+ \Hess_{\rr_\pp}(u_2,u_2)\le 0$$
for $u\in \cup_{x\in\pp M}T_xM$.
Since $\Hess_{\rr_\pp }$ is smooth on the compact set $\pp_0 M$, this implies
$$\Hess_{\rr_\pp}(u,u)\le c|u|^2\rr_\pp(x),\ \ x\in M, u\in T_xM$$
for some  constant $c>0$. Then the desired estimate \eqref{HESS} follows from
$$\Hess_{\log \rr_\pp}(u,u)= \rr_\pp^{-1}\Hess_{\rr_\pp}(u,u)-\rr_\pp^{-2}\<\nn\rr_\pp, u\>^2\le c |u|^2,\ \ u\in TM^\circ.$$
\end{proof}

 By Lemma \ref{L3.0}, when $\pp M$ is convex, there exists a constant $K\ge 0$ such that
\beq\label{CV} \Ric-\Hess_{V+2\log\phi_0}\ge -K.\end{equation}
 Since the diffusion process generated by $L_0:=\DD+\nn(V+2\log\phi_0)$ is non-explosive in $M^\circ$,  by \eqref{CV} and  Bakry-Emery's semigroup calculus,
 (see for instance \cite{BGL} or \cite[Theorem 2.3.3]{W14}), we have
 \beq\label{SM1} |\nn P_t^0g|\le \e^{Kt} P_t^0|\nn g|,\ \ t\ge 0, g\in C_b^1(M)\end{equation}  and for any $p>1$, there exists a constant $c(p)>0$ such that
 \beq\label{SM2}\beg{split}  |\nn P_t^0g|^2&\le \ff{2K \{(P_t^0|g|^{p\land 2} )(P_t^0|g|)^{(2-p)^+} -(P_t^0|g|)^{2}\} }{(p\land 2)(p\land 2-1) (1-\e^{-2Kt})}\\
 &\le \ff{c(p)}{1\land t} (P_t^0|g|^p)^{\ff 2 p},\ \ t>0, g\in \B_b(M).\end{split}\end{equation}

When $\pp M$ is non-convex, we take as in \cite{W07}    a conformal change of metric to make it convex under the new metric. More precisely, we have the following result.

 \beg{lem}\label{CCV} There exists a function $1\le \phi\in C_b^\infty(M)$ such that  $\pp M$ is convex under the metric $\<\cdot,\cdot\>_\phi:=\phi^{-2}\<\cdot,\cdot\>.$ Moreover,
 there exists a smooth vector field $Z_\phi$ on $M$ such that
\beq\label{LLL} L_0= \phi^{-2}\DD^\phi+Z_\phi+ 2\phi^{-1}\nn^\phi\log\phi_0,\end{equation}
where $\nn^\phi$ and $\DD^\phi$ are the  gradient  and Lapalce-Beltrami  operators induced by $\<\cdot,\cdot\>_\phi$ respectively.\end{lem}

\beg{proof}   let $\dd>0$ such that the second fundamental form of $\pp M$ is bounded below by $-\dd$. Take $1\le \phi\in C_b^\infty(M)$ such that
$\phi= 1+\dd \rr_\pp$  in a neighborhood of $\pp M$ in which the distance function $\rr_\pp$ to $\pp M$ is smooth.
By \cite[Lemma 2.1]{W14b}(see also \cite{W07}),   $\pp M$ is convex under the metric $\<\cdot,\cdot\>_\phi:=\phi^{-2}\<\cdot,\cdot\>.$
Next, according to  the proof of \cite[Lemma 2.2]{W14b},
there exists a smooth vector field $Z_\phi$ on $M$ such that \eqref{LLL} holds.
\end{proof}

 Let $1\le \phi\in C_b^\infty(M)$ be as in Lemma \ref{CCV}, and let $P_t^\phi$ be the diffusion semigroup generated by
$$L^\phi:= \phi L_0= \phi^{-1} \DD^\phi+ \phi Z_\phi + 2 \nn^\phi\log\phi_0.$$
We have the following result.

\beg{lem}\label{GRD}  Let $1\le \phi\in C_b^\infty(M)$ be as in Lemma \ref{CCV}.
\beg{enumerate} \item[$(1)$]  For any $q\in (1,\infty]$,  there exists a constant $c(q) >0$ such that
\beq\label{1*3} |\nn^\phi P_t^\phi f|_\phi\le \ff {c(q)} {\ss t}  (P_t^\phi  |f|^q)^{\ff 1 q},\ \ t>0, f\in C_b^1(M).\end{equation} Moreover,  there exists a constant $K>0$ such that
\beq\label{1*2} |\nn^\phi P_t^\phi f|_\phi\le \e^{K t} P_t^\phi|\nn^\phi f|_\phi,\ \ t>0, f\in C_b^1(M).\end{equation}
\item[$(2)$] There exists a constant $c>0$ such that
\beq\label{PHI} \|P_t^\phi\|_{L^p(\mu_0)\to L^\infty(\mu_0)}\le c (1\land t)^{-\ff{d+2}{2p}},\ \ t>0,p\in [1,\infty].\end{equation}\end{enumerate} \end{lem}

\beg{proof}  (1) Since $\pp M$ is convex under the metric $\<\cdot,\cdot\>_\phi$, by Lemma \ref{L3.0}, we find a constant $K_0^\phi>0$ such that
\beq\label{HS} 2\Hess^\phi_{\log\phi_0} (u,u)\le K_0^\phi|u|_\phi^2,\ \ u\in TM^\circ,\end{equation}
where $\Hess^\phi$ is the Hessian tensor induced by the metric $\<\cdot,\cdot\>_\phi$.
Since the operator $A^\phi:=\phi^{-1} \DD^\phi +\phi Z_\phi$ is a $C^2$-smooth strictly elliptic second order differential operator on the compact manifold $M$,
it has bounded below Bakry-Emery curvature; that is,  there exists a constant $K_1^\phi>0$ such that
$$A^\phi |\nn^\phi f|_\phi^2 - 2 \<\nn^\phi A^\phi f,\nn^\phi f\>_\phi \ge - K_1^\phi |\nn^\phi f|_\phi^2,\ \ f\in C^\infty(M).$$
Combining this with   \eqref{HS} we obtain
$$ L^\phi |\nn^\phi f|^2_\phi - 2 \<\nn^\phi L^\phi f,\nn^\phi f\>_\phi
\ge -(K_0^\phi+ K_1^\phi)|\nn^\phi f|^2_\phi=:-K^\phi |\nn^\phi f|^2_\phi,\ \ f\in C^\infty(M^\circ),$$ which  means that the Bakry-Emery curvature of $L^\phi$ is bounded below by $-K^\phi$.
By the same reason leading to
 \eqref{SM1} and \eqref{SM2}, this implies \eqref{1*3} and \eqref{1*2}.

(2) To estimate  $\|P_{t}^\phi\|_{L^p(\mu_0)\to L^\infty(\mu_0)},$ we make use of \cite[Theorem 4.5(b)]{W00} or \cite[Theorem 3.3.15(2)]{Book}, which says that \eqref{PRR0} implies the super Poincar\'e inequality
$$\mu_0(f^2)\le r \mu_0(|\nn f|^2)+\bb (1+ r^{-\ff{d+2}2})\mu_0(|f|)^2,\ \ r>0, f\in C_b^1(M)$$ for some constant $\bb>0$.
Let $\mu^\phi= \ff{\phi^{-1}\mu_0}{\mu_0(\phi^{-1})}$. By $L^\phi=\phi L_0$ we obtain
$$\EE^\phi(f,g):= -\int_M f L^\phi g \d\mu^\phi= -\ff 1 {\mu_0(\phi^{-1})} \int_M fL_0 g\d\mu_0= \ff 1 {\mu(\phi^{-1})} \mu_0(\<\nn f,\nn g\>),
\ \ f,g\in C_b^2(M).$$
Then the above super Poincar\'e inequality implies
$$\mu^\phi (f^2)\le r \EE^\phi(f,f) +\bb' (1+ r^{-\ff{d+2}2})\mu^\phi(|f|)^2,\ \ f\in C_b^1(M)$$
for some constant $\bb'>0$. Using \cite[Theorem 4.5(b)]{W00} or \cite[Theorem 3.3.15(2)]{Book} again, this implies
$$\|P_t^\phi\|_{L^p(\mu^\phi)\to L^\infty(\mu^\phi)}\le \kk (1\land t)^{-\ff{d+2}{2p}},\ \ t>0, p\in [1,\infty]$$ for some constant $\kk>0$. Noting that
$$\|\phi\|_\infty^{-1} \mu_0\le \mu^\phi\le \|\phi\|_\infty \mu_0,$$ we find a constant $c>0$ such that \eqref{PHI} holds.
\end{proof}

\beg{lem}\label{LFW} For any $p\in (1,\infty]$ and $q\in (1,p)$, there exists a constant $c>0$ such that for any $f\in D(L_0)$,
\beq\label{KNN} \|\nn P_t^0f\|_\infty\le c \e^{-\ll_0 t} \big\{(1\land t)^{-\ff 1 2 -\ff{d+2}{2p}} \|f\|_{L^p(\mu_0)}+ (1\land t)^{\ff 1 2 - \ff{q(d+2)}{2p}} \|L_0f\|_{L^{p}(\mu_0)}\big\},\ \ t>0.\end{equation}
Consequently, there exists a constant $c>0$ such that
\beq\label{CC'} \|\nn (\phi_m\phi_0^{-1})\|_\infty\le c m^{\ff{d+4}{2d}},\ \ m\ge 1.\end{equation}\end{lem}

\beg{proof}  (a) By the semigroup property and the $L^p(\mu_0)$ contraction  of $P_t^0$,     for the proof of  \eqref{KNN} it suffices to consider  $t\in (0,1].$
  Since   $1\le \phi\in C_b^\infty(M)$, we have $\D(L_0)=\D(L^\phi)$ and
\beq\label{FW1} P_t^0f= P_t^\phi f - \int_0^t P_s^\phi\{(\phi-1)P_{t-s}^0L_0 f\}\d s,\ \ t\ge 0, f\in \D(L_0).\end{equation}
Next, by \eqref{1*3} and \eqref{PHI}, we obtain 
\beq\label{FW2}\beg{split} &\|\nn P_t^\phi f\|_\infty = \|\nn P_{t/2}^{\phi}(P_{t/2}^\phi f)\|_\infty\\
&\sm  t^{-\ff 1 2} \|P_{t/2}^\phi f\|_\infty\sm  t^{-\ff 1 2-\ff{d+2}{2p}}\|f\|_{L^p(\mu_0)},\ \ t\in (0,1].\end{split} \end{equation}
Combining this with \eqref{**1} and \eqref{1*3} leads to 
\beg{align*} &\int_0^t \|\nn P_s^\phi \{(\phi-1) P_{t-s}^0 L_0f\}\|_\infty \d s\sm  \int_0^t s^{-\ff 1 2} \big\|\{P_s^\phi |P_{t-s}^0 L_0f|^{q}\}^{\ff 1 q}\big\|_\infty\d s\\
&\sm \int_0^{\ff t 2} s^{-\ff 1 2} \|P_{t-s}^0L_0f\|_\infty\d s+     \int_{\ff t 2}^t s^{-\ff 1 2}
\big\|\{P_s^\phi |P_{t-s}^0L_0f|^{q}\}^{\ff 1 q}\big\|_\infty\d s\\
&\sm  \int_0^{\ff t 2}s^{-\ff 1 2}  \|P_{t-s}^0\|_{L^p(\mu_0)\to L^\infty (\mu_0)}\|L_0f\|_{L^p(\mu_0)}  \d s
+  \int_{\ff t 2}^t s^{-\ff 1 2} \|P_s^\phi\|_{L^{\ff p q}(\mu_0)\to L^\infty (\mu_0)} \|L_0f\|_{L^{p}(\mu_0)} \\
  &\sm  t^{\ff 1 2 -\ff{q(d+2)}{2p}}\|L_0f\|_{L^{p}(\mu_0)}.\end{align*}
Substituting this and  \eqref{FW2} into   \eqref{FW1},  we prove \eqref{KNN}.

(b) Applying \eqref{KNN} to $p=\infty$, $f=\phi_m\phi_0^{-1}, t= (\ll_m-\ll_0)^{-1}$ and using  \eqref{PR1},  we obtain
$$\e^{-1} \|\nn(\phi_m\phi_0^{-1})\|_\infty  \sm (\ll_m-\ll_0)^{\ff 1 2} \|\phi_m\phi_0^{-1}\|_\infty,\ \ m\ge 1. $$
 This together with \eqref{EG} and \eqref{CC} implies \eqref{CC'} for some constant $c>0$.  \end{proof}

 \section{Upper bound estimate}

 According to \cite[Lemma 2.3]{WZ20}, we have
 \beq\label{2.2} \W_2(\mu_t^\nu,\mu_0)^2\le \int_M \ff{ |\nn L_0^{-1}(h_t^\nu-1)|^2}{\scr M(h_t^\nu,1)}\d\mu_0a                                                                                                                                 ,\end{equation}
 where
 $$h_t^\nu:=\ff{\d\mu_t^\nu}{\d\mu_0},\ \  \scr M (a,b):= 1_{\{a\land b>0\}}\ff{a-b}{\log a-\log b}.$$
 So, to investigate the upper bound estimate, we first calculate $h_t^\nu$.

By \eqref{ON*}, we have   
\beq\label{PSIS}\psi_s^\nu:= \int_M\phi_0(x) p_s^0(x,\cdot)\nu(\d x)=\nu(\phi_0)+\sum_{m=1}^\infty \nu(\phi_m) \e^{-(\ll_m-\ll_0)s} \phi_m\phi_0^{-1},\ \ s>0.\end{equation}
Next,  \eqref{PR0} and \eqref{ON*} imply
 \beq\label{NB*} \nu(P_s^D f) =\e^{-\ll_0 s} \nu(\phi_0P_s^0(f\phi_0^{-1})) =\e^{-\ll_0 s} \int_M \psi_s^\nu\phi_0^{-1} f\d \mu_0,\ \ f\in \B^+(M),\end{equation}
 where $\B^+(M)$ is the class of nonnegative measurable functions on $M$.
Moreover, for any $t\ge s> 0$, by the Markov property, \eqref{2.1}, \eqref{PR0} and \eqref{NB*},  we obtain
 \beg{align*} & \int_M f \d \E^\nu[\dd_{X_s}1_{\{t<\tau\}}]= \E^\nu\big[f(X_s) 1_{\{s<\tau\}} (P_{t-s}^D 1)(X_s)\big]= \nu \big(P_s^D\{fP_{t-s}^D1\} \big)\\
 &=  \e^{-\ll_0t} \int_M( \psi_s^v P_{t-s}^0 \phi_0^{-1}) f \d\mu_0,\ \  f\in \B^+(M).\end{align*}
 Then
 $$\ff{\d \E^\nu[\dd_{X_s}1_{\{t<\tau\}}]}{\d\mu_0}= \e^{-\ll_0t}  \psi_s^v P_{t-s}^0\phi_0^{-1}.$$
 Noting that \eqref{NB*} implies
 $$\E^\nu [1_{\{t<\tau\}}]= \nu(P_t^D1)= \e^{-\ll_0 t}\mu_0(\psi_t^\nu \phi_0^{-1})= \e^{-\ll_0 t}\nu(\phi_0 P_t^0\phi_0^{-1}),$$
we arrive at
 \beq\label{*10}   \beg{split} & h_t^\nu:= \ff{\d\mu_t^\nu}{\d\mu_0}= \ff {1} {t   \E^\nu 1_{\{t<\tau\}}}\int_0^t \ff{\d \E^\nu [\dd_{X_s}1_{\{t<\tau\}}]}{\d\mu_0} \d s=
 1+\rr_t^\nu,  \\
& \rr_t^\nu:=  \ff 1 {t\nu(\phi_0 P_t^0\phi_0^{-1})} \int_0^t \big\{\psi_s^\nu P_{t-s}^0\phi_0^{-1} -\nu(\phi_0P_t^0\phi_0^{-1}) \big\}\d s.\end{split}\end{equation}
By \eqref{**1}, $\|\phi_0\|_\infty<\infty$   and $\|\phi_0^{-1}\|_{L^2(\mu_0)}=1$, we find a constant $c>0$ such that
\beq\label{GY3}\beg{split}& |\nu(\phi_0 P_t^0\phi_0^{-1})- \nu(\phi_0)\mu(\phi_0)| \le \nu(\phi_0)\|P_t^0\phi_0^{-1}-\mu_0(\phi_0^{-1})\|_\infty\\
&\le c\e^{-(\ll_1-\ll_0)t},\ \ t\ge 1,\nu\in \scr P_0.\end{split} \end{equation}

Due to the lack of simple  representation of the product $\psi_s^\nu P_{t-s}^0\phi_0^{-1}$ in terms of the eigenbasis $\{\phi_m\phi_0^{-1}\}_{m\ge 0}$, it  is inconvenient to estimate   the upper bound in \eqref{2.2}.
To this end, below we  reduce this product  to a linear combination of  $\psi_s^\nu $ and $P_{t-s}^0\phi_0^{-1}$,
for which the spectral representation works.
Write
\beq\label{GY1} \beg{split} & \psi_s^\nu P_{t-s}^0 \phi_0^{-1} -\nu(\phi_0P_t^0\phi_0^{-1})= I_1(s) + I_2(s),\\
&I_1(s):= \{\psi_s^\nu-\nu(\phi_0)\}\cdot\{P_{t-s}^0 \phi_0^{-1} -\mu(\phi_0)\} +\nu(\phi_0\{\mu(\phi_0)-P_t^0\phi_0^{-1}\}),\\
&I_2(s):= \mu(\phi_0) \{\psi_s^\nu-\nu(\phi_0)\}+\nu(\phi_0) \{ P_{t-s}^0 \phi_0^{-1} -\mu(\phi_0)\}.\end{split}\end{equation}
By   \eqref{ONN},   \eqref{ON*} and \eqref{PSIS},  we have
\beq\label{GYY} \beg{split}  & P_{t-s}^0 \phi_0^{-1} -\mu(\phi_0)=  \sum_{m=1}^\infty \mu(\phi_m) \e^{-(\ll_m-\ll_0)(t-s) } \phi_m\phi_0^{-1},\\
& \psi_s^\nu-\nu(\phi_0)  =   \sum_{m=1}^\infty \nu(\phi_m) \e^{-(\ll_m-\ll_0)s } \phi_m\phi_0^{-1},\ \ t>s>0. \end{split}\end{equation}
Then
\beq\label{GY2} \beg{split} &\rr_t^\nu = \tt\rr_t^\nu +  \ff 1 {t\nu(\phi_0 P_t^0\phi_0^{-1})} \int_0^t I_1(s)\d s-A_t,\\
&\tt\rr_t^\nu:= \ff 1 {t\nu(\phi_0 P_t^0\phi_0^{-1})} \sum_{m=1}^\infty\ff{ \mu(\phi_0) \nu(\phi_m)+\nu(\phi_0)\mu(\phi_m)}{\ll_m-\ll_0} \phi_m\phi_0^{-1}\\
&A_t:= \ff 1 {t\nu(\phi_0 P_t^0\phi_0^{-1})}\sum_{m=1}^\infty\ff{\{ \mu(\phi_0) \nu(\phi_m)+\nu(\phi_0)\mu(\phi_m)\}\e^{-(\ll_m-\ll_0)t} }{\ll_m-\ll_0} \phi_m\phi_0^{-1}.\end{split}\end{equation}
Since $\rr_t^\nu\in L^1(\mu_0)$, the following lemma implies  $\tt\rr_t^\nu\in L^1(\mu_0)$ for $t>0.$

\beg{lem}\label{LGY1}  For any $t_0>0$, there exists a constant $c>0$ such that
\beq\label{GY4} \mu_0(|\rr_t^\nu-\tt\rr_t^\nu|)\le c\|h\|_{L^2(\mu)} \e^{-(\ll_1-\ll_0)t},\ \ t\ge t_0, \nu=h\mu\in\scr P_0.\end{equation}\end{lem}
\beg{proof} By \eqref{EG} and \eqref{CC}, for any $t_0>0$  we have 
\beq\label{AT} \sum_{m=1}^\infty\|\phi_m\|_\infty\e^{-(\ll_m-\ll_0)t}\sm  \e^{-(\ll_1-\ll_0)t},\ \ t\ge t_0.\end{equation} Combining this with \eqref{GY2} and \eqref{GY3}, and noting that
$\|h\phi_0^{-1}\|_{L^2(\mu_0)}= \|h\|_{L^2(\mu)}, $  it suffices to show that
 \beq\label{AUU}B:= \ff 1 t\int_0^t \big\| \{\psi_s^\nu-\nu(\phi_0)\}\cdot\{P_{t-s}^0\phi_0^{-1}-\mu(\phi_0)\}\big\|_{L^1(\mu_0)}\d s\sm \|h\|_{L^2(\mu)}\e^{-(\ll_1-\ll_0)t},\ \ t\ge t_0.\end{equation}  
Since $\|\phi_0^{-1}\|_{L^2(\mu_0)}=1$ and $\psi_s^\nu= P_0^s(h\phi_0^{-1})$ for $\nu= h\mu$,   \eqref{000} yields  that 
 \beg{align*}  
 &B\le  \ff {1} {t} \int_0^t   \|P_{t-s}^0\phi_0^{-1}- \mu_0(\phi_0^{-1})\|_{L^2(\mu_0)} \|P_{s}^0(h\phi_0^{-1})- \mu_0(h\phi_0^{-1})\|_{L^{2}(\mu_0)} \d s\\
&\le \ff {1}t \int_0^t  \|P_{t-s}^0-\mu_0\|_{L^{2}(\mu_0)} \|P_{s}^0-\mu_0\|_{L^2(\mu_0)}\|h\|_{L^2(\mu)} \d s\\
&\sm \|h\|_{L^2(\mu)}\e^{-(\ll_1-\ll_0)t},\ \ t\ge t_0.\end{align*}
\end{proof}

\beg{lem}\label{LGY2}  For any    $\aa>0$, there exist    constants $c_0,t_0 >0$ such that
\beq\label{GY4} \tt\rr_{t}^\nu\ge -\ff{c_0} {\nu(\phi_0)t},\ \ t\ge t_0,\ \ \nu\in \scr P_0, \nu\in \scr P_0.  \end{equation}
Consequently,  if $\nu=h\mu$ with $h\in L^2(\mu)$, then $\tt\mu_t^\nu:=(1+\tt\rr_t^\nu)\mu_0$ is a probability measure for  $t>t_0(1+c_0)$.\end{lem}
\beg{proof}
By Lemma \ref{LGY1},  if $\nu=h\mu$ with $h\in L^2(\mu)$, we have  $\tt\rr_t^\nu\in L^1(\mu_0)$ for  $t>0$, and it is easy to see that $\mu_0(\tt\rr_t^\nu)=0$.
Since \eqref{GY4} implies $1+\tt\rr_t^\nu>0$ for   $t>t_0(1+c_0)$,     $\tt\mu_t^\nu$ is a probability measure.
It remains to prove  \eqref{GY4}.

By \eqref{GY3} and \eqref{GY2}, it suffices to find a constant $c_1>0$ such that
\beq\label{*BB} g:=\sum_{m=1}^\infty\ff{ \mu(\phi_0) \nu(\phi_m)+\nu(\phi_0)\mu(\phi_m)}{\ll_m-\ll_0} \phi_m\phi_0^{-1}\ge -c_1.\end{equation}
By \eqref{EG} and \eqref{CC}, we have 
\beq\label{*ABK} \|P_1^0g\|_\infty \le c_2:= \sum_{m=1}^\infty\ff{2\|\phi_0\|_\infty \|\phi_m\|_\infty \|\phi_m\phi_0^{-1}\|_\infty}{(\ll_m-\ll_0)\e^{\ll_m-\ll_0}}<\infty.\end{equation}
Next, by \eqref{GYY} and the same formula for $\mu=\nu$, we obtain
\beq\label{999} P_s^0g= (-L_0)^{-1} \big\{\mu(\phi_0)(\psi_s^\nu-\nu(\phi_0))+\nu(\phi_0)(\psi_s^\mu-\mu(\phi_0))\big\} = (-L_0)^{-1} g_s,\ \ s>0,\end{equation}
where by $\phi_0,\psi_s^\nu,\psi_s^\mu\ge 0$,
$$g_s:= \mu(\phi_0)(\psi_s^\nu-\nu(\phi_0))+\nu(\phi_0)(\psi_s^\mu-\mu(\phi_0))\ge -2 \mu(\phi_0)\nu(\phi_0)\ge -2\nu(\phi_0),\ \ s>0.$$
This together with \eqref{999} yields
$$-L_0 P_s^0 g\ge -2\nu(\phi_0),\ \ s>0.$$
Therefore, it follows from \eqref{*ABK} that
$$g= P_1^0 g- \int_0^1L_0 P_r^0g\d r \ge -c_2 -2\nu(\phi_0)\ge -c_2-2\|\phi_0\|_\infty.$$
So,  \eqref{*BB} holds for  $c_1= c_2+2\|\phi_0\|_\infty$. \end{proof}

\beg{lem}\label{LAY4}   There exist    constants $c,t_0>0$ such that  for any $t\ge t_0$,  and any $\nu\in \scr P_0$ with $\nu=h\mu$ such that $h\in L^2(\mu)$, 
we have $\tt\mu_t^\nu \in\scr P_0$ and 
\beq\label{KKL}  t^2\W_2(\tt\mu_{t}^\nu, \mu_0)^2 \le  \ff{1+ct^{-1}} {\{\mu(\phi_0)\nu(\phi_0)\}^2} \sum_{m=1}^\infty \ff{\{\nu(\phi_0)\mu_0(\phi_m)+ \mu(\phi_0) \nu(\phi_m)\}^2}{(\ll_m-\ll_0)^3}.\end{equation}
\end{lem}
\beg{proof}     By Lemma \ref{LGY2}, there exist constants  $c,t_0>0$ such that $\tt\mu_t^\nu \in\scr P_0$ for $t\ge t_0$, and
$$ \scr M(1+\tt\rr_t^\nu,1)\ge 1\land (1+\tt\rr_t^\nu) \ge \ff 1 {1+ct^{-1}},\ \ t\ge t_0.$$
So,   \cite[Lemma 2.3]{WZ20} implies
 \beq\label{2.22}\beg{split} &  \W_2(\tt\mu_{t}^\nu, \mu_0)^2 \le    \int_M \ff{ |\nn L_0^{-1}\tt \rr_{t}^\nu|^2}{\scr M(1+\tt\rr_{t}^\nu,1)}\d\mu_0  \le
 (1+ct^{-1})  \mu_0(|\nn L_0^{-1}\tt \rr_{t}^\nu|^2),\ \  t\ge t_0.\end{split} \end{equation}
Next,   \eqref{PR1} and \eqref{GY2} yield
$$ t^2\mu_0(|\nn L_0^{-1}\tt \rr_{t}^\nu|^2)= \ff 1 {\{\nu(\phi_0 P_t^0\phi_0^{-1})\}^2}\sum_{m=1}^\infty \ff{ \{\mu(\phi_0)\nu(\phi_m) + \nu(\phi_0)\mu(\phi_m)\}^2}{(\ll_m-\ll_0)^3 }.$$
Combining this with \eqref{GY3} and \eqref{2.22}, we finish the proof.
\end{proof}

We are now ready to prove the following result.

\beg{prp}\label{PA} For any $\nu\in\scr P_0$, 
\beq\label{LOO} \limsup_{t\to\infty}  \big\{t^2\W_2(\mu_t^\nu,\mu_0)^2\big\}\le I.\end{equation} \end{prp}

 \beg{proof}   
(1)  We first consider $\nu=h\mu$ with $h\in L^2(\mu)$.  Let $D$ be the diameter of $M$. By Lemma \ref{LGY1}, there exist  constants $c_1,t_0>0$ such that
$\tt\mu_t^\nu$ is probability measure for $t\ge t_0$ and 
 \beq\label{OK} \W_2(\mu_t^\nu,\tt\mu_t^\nu)^2 \le D^2 \|\mu_t^\nu-\tt\mu_t^\nu\|_{var}=D^2\mu_0(|\rr_t^\nu-\tt\rr_t^\nu|) \le c_1 \|h\|_{L^2(\mu)}\e^{-(\ll_1-\ll_0)t},\ \ t\ge t_0.\end{equation}
Combining this with Lemma \ref{LAY4} 
and the   triangle inequality of $\W_2$, 
we obtain
\beq\label{*AG} t^2\W_2(\mu_t^\nu, \mu_0)^2 \le (1+\dd^{-1}) c_1 t^2 \e^{-(\ll_1-\ll_0)t} \|h\|_{L^2(\mu)} +(1+\dd) (1+ct^{-1}) I,\ \ \dd>0.\end{equation}

(2) In general, we may go back to the first situation by shifting a small time $\vv>0$. More precisely, by the Markov property, \eqref{2.1},  \eqref{PR0} and \eqref{PSIS},
for any $f\in \B_b(M)$ and $t\ge s\ge \vv>0,$ we have 
\beg{align*} &\E^\nu[f(X_s) 1_{\{t<\tau\}}] = \E^\nu\big[1_{\{\vv<\tau\}} \E^{X_\vv} (f(X_{s-\vv})1_{\{t-\vv<\tau\}})\big] \\
&=\int_{M} p_\vv^D(x,y) \E^y[f(X_{s-\vv})1_{\{t-\vv<\tau\}}]\nu(\d x)\mu(\d y) \\
&= \e^{-\ll_0\vv} \int_M (\psi_\vv^\nu\phi_0)(y) \E^y[f(X_{s-\vv})1_{\{t-\vv<\tau\}}]\nu(\d x)\mu(\d y).\end{align*}
With $f=1$ this implies 
$$\P^\nu(t<\tau)=  \e^{-\ll_0\vv} \int_M (\psi_\vv^\nu\phi_0)(y) \P^y(t-\vv<\tau)\mu(\d y)\mu(\d y).$$
So, letting 
$$\nu_\vv=\ff{\psi_\vv^\nu\phi_0}{\mu(\psi_\vv^\nu\phi_0)}=:h_\vv\mu,$$
we arrive at 
$$\E^\nu[f(X_s)|t<\tau] =\ff{\E^\nu[f(X_s) 1_{\{t<\tau\}}] }{\P^\nu(t<\tau)} = \ff{\E^{\nu_\vv}[f(X_{s-\vv}) 1_{\{t-\vv<\tau\}}] }{\P^{\nu_\vv}(t-\vv<\tau)}=
\E^{\nu_\vv}[f(X_{s-\vv})|t-\vv<\tau].$$
Therefore, 
\beq\label{*AGG} \mu_{t,\vv}^\nu:= \ff 1 {t-\vv}\int_\vv^t \E^\nu(\dd_{X_s}|t<\tau)\d s= \mu_{t-\vv}^{\nu_\vv},\ \ t>\vv.\end{equation}
Since 
$$\mu(\psi_\vv^\nu\phi_0)= \int_M p_\vv^0(x,y) \phi_0(x)\phi_0(y) \nu(\d x) \mu(\d y) =\nu(\phi_0 P_\vv^0\phi_0^{-1})\ge \nu(\phi_0)\|\phi_0\|_\infty^{-1}=:\aa>0,$$
by \eqref{PRR0} we find a constant $c_2>0$ such that 
\beq\label{MX0} \|h_\vv\phi_0^{-1}\|_{L^2(\mu_0)} \le \aa^{-1} \|\psi_\vv^\nu\|_{L^2(\mu_0)} \le  \aa^{-1} \|\phi_0\|_\infty   \|p_{\vv}^0\|_{L^\infty(\mu_0)} \le c_2  \vv^{-\ff{d+2}2},\ \ \vv\in (0,1).\end{equation}
Then \eqref{*AG} and \eqref{*AGG} yield
\beq\label{*AG2} \beg{split} &t^2\W_2(\mu_{t,\vv}^\nu,\mu_0)^2 \\
&\le (1+\dd^{-1}) c_1c_2\aa^{-1}  t^2 \e^{-(\ll_1-\ll_0)t} \vv^{-\ff{d+2}2} +(1+\dd) (1+ct^{-1}) I_\vv,\ \ \dd>0, \vv\in (0,1),\end{split}\end{equation}
where 
$$I_\vv:= \ff1 {\{\mu(\phi_0)\nu_\vv(\phi_0)\}^2} \sum_{m=1}^\infty \ff{\{\nu_\vv(\phi_0)\mu(\phi_m)+ \mu(\phi_0) \nu_\vv(\phi_m)\}^2}{(\ll_m-\ll_0)^3}.$$
By \eqref{PR0}, \eqref{PR1} and \eqref{PSIS}, we have  
\beg{align*} &\mu(\psi_\vv^\nu \phi_0) =\nu(\phi_0P_\vv^{-1}\phi_0^{-1})=\e^{\ll_0\vv} \nu(P_\vv^D1),\\
&\mu(\psi_\nu\phi_0)= \nu(\phi_0P_\vv^0(\phi_m\phi_0^{-1})) =\e^{-(\ll_m-\ll_0)\vv} \nu(\phi_m),\end{align*}
so that 
 $$  \nu_\vv(\phi_m)= \ff{\e^{-\ll_m \vv} \nu(\phi_m)}{\nu(P_\vv^D1)},\ \ m\ge 0. $$
Thus, $\lim_{\vv\to 0} \nu_\vv(\phi_0)=\nu(\phi_0)$ and there exists a constant $C>1$ such that
\beq\label{BNN} C^{-1} \e^{-\ll_m \vv} |\nu(\phi_m)| \le |\nu_\vv(\phi_m)|\le C|\nu(\phi_m)|,\ \ m\ge 1,\vv\in (0,1).\end{equation} 
Therefore, if $I<\infty$, by this and 
\beq\label{DMM} \sum_{m=1}^\infty \mu(\phi_m)^2\le \mu(1)=1,\end{equation}  we may apply the dominated convergence theorem to derive  $\lim_{\vv\to 0} I_\vv=I$. On the other hand, if 
$I=\infty$, which is equivalent to 
$$ \sum_{m=1}^\infty \ff{\nu(\phi_m)^2}{(\ll_m-\ll_0)^3}=\infty,$$
 then by \eqref{BNN} and  the monotone convergence theorem we get
$$\liminf_{\vv\to 0} \sum_{m=1}^\infty \ff{\nu_\vv(\phi_m)^2}{(\ll_m-\ll_0)^3}\ge C^{-2} \liminf_{\vv\to 0} \sum_{m=1}^\infty \ff{\e^{-2\ll_m\vv}\nu(\phi_m)^2}{(\ll_m-\ll_0)^3}=\infty,$$ which together with \eqref{DMM} and $\nu_\vv(\phi_0)\to\nu(\phi_0)$ implies
\beg{align*} &\liminf_{\vv\to 0} I_\vv =\ff1 {\{\mu(\phi_0)\nu(\phi_0)\}^2}  \liminf_{\vv\to 0}  \sum_{m=1}^\infty \ff{\{\nu_\vv(\phi_0)\mu(\phi_m)+ \mu(\phi_0) \nu_\vv(\phi_m)\}^2}{(\ll_m-\ll_0)^3}\\
&\ge \ff1 {\{\mu(\phi_0)\nu(\phi_0)\}^2}  \liminf_{\vv\to 0} \ff{ \ff 1 2 \{\mu(\phi_0) \nu_\vv(\phi_m)\}^2-\|\phi_0\|_\infty^2 \mu(\phi_m)^2}{(\ll_m-\ll_0)^3}=\infty.\end{align*}
In conclusion, we have 
\beq\label{LMT} \lim_{\vv\to 0}I_\vv=I.\end{equation}This together with \eqref{*AG2} for $\vv= t^{-2}$ gives
\beq\label{*AGN} \limsup_{t\to\infty} \big\{t^2\W_2(\mu_{t,t^{-2}}^\nu,\mu_0)^2\big\}\le I.\end{equation}
On the other hand, it is easy to see that
$$\|\mu_{t,\vv}^\nu-\mu_t^\nu\|_{var}\le \ff{2\vv} t,\ \ 0<\vv<t,$$ so that 
\beq\label{*AG3} \W_2(\mu_t^\nu,\mu_{t,t^{-2}}^\nu)^2\le D^2 \|\mu_{t,t^{-2}}^\nu-\mu_t^\nu\|_{var}\le 2D^2 t^{-3},\ \ t>1.\end{equation}
Combining this with \eqref{*AGN}, we prove \eqref{LOO}.  
  \end{proof}

\section{Lower bound estimate and the finiteness of the limit}

We will  follow the idea of \cite{AMB, WZ20}, for which we need to modify $\tt\mu_t^\nu$ as follows. For any $\bb>0$, consider
$$\tt\mu^\nu_{t,\bb}:= (1+\tt\rr^\nu_{t,\bb})\mu_0,\ \ \tt\rr^\nu_{t,\bb}:= P_{t^{-\bb}}^0 \tt\rr_t^\nu,\ \ t>0.$$
According to Lemma \ref{LGY2}, there exists $t_0>0$ such that
\beq\label{ALS}  \tt h_t^\nu :=1+\tt \rr_t^\nu\ge \ff 1 2,\ \ \tt h_{t,\bb}^\nu:=1+\tt\rr_{t,\bb}^\nu\ge \ff 12,\ \ \bb>0, t\ge t_0. \end{equation}
Consequently,  $\tt\mu_{t,\bb}^\nu$ and $\tt\mu_t^\nu$ are probability measures for any $\bb>0, t\ge t_0.$

\beg{lem}\label{L3.1}  For any $\bb>0$, there exists a  constant $c>0$ such that  $f_{t,\bb}:= L_0^{-1}\tt\rr_{t, \bb}^\nu$   satisfies
$$\|f_{t,\bb}\|_\infty+\|L_0f_{t,\bb}\|_\infty+\|\nn f_{t,\bb}\|_\infty\le c t^{\ff{5\bb d} 4-1},\ \ t\ge 1. $$
\end{lem}

\beg{proof} By \eqref{PR1} and \eqref{GY2}, we have
\beg{align*} &f_{t,\bb}= -   \sum_{m=1}^\infty\ff{ \{\mu(\phi_0)\nu(\phi_m) + \nu(\phi_0)\mu(\phi_m)\} \e^{-(\ll_m-\ll_0)t^{-\bb}}} {t(\ll_m-\ll_0)^2\nu(\phi_0P_t^0\phi_0^{-1})}  \big(\phi_m\phi_0^{-1}\big),\\
& L_0 f_{t,\bb}= \sum_{m=1}^\infty\ff{ \{\mu(\phi_0)\nu(\phi_m) + \nu(\phi_0)\mu(\phi_m)\} \e^{-(\ll_m-\ll_0)t^{-\bb}}} {t(\ll_m-\ll_0)\nu(\phi_0P_t^0\phi_0^{-1})}  \big(\phi_m\phi_0^{-1}\big).\end{align*}
 Combining these with \eqref{EG},  \eqref{CC},    \eqref{GY3},   and
 $$|\mu(\phi_0)\nu(\phi_m) + \nu(\phi_0)\mu(\phi_m)|\le \|\phi_0\|_\infty +\|\phi_m\|_\infty\sm m,\ \ m\ge 1, $$   
  we find  a constant  $c_1 >0$  such that
 \beg{align*} &t\{\|f_{t,\bb}\|_\infty+\|L_0f_{t,\bb}\|_\infty\}\sm  \sum_{m=1}^\infty\ff{  \e^{-(\ll_m-\ll_0)t^{-\bb}} m^{\ff{3d+2}{2d} }}{\ll_m-\ll_0 }   \\
 &\sm   \sum_{m=1}^\infty   \e^{-c_1m^{\ff 2 d} t^{-\bb} } m^{\ff{3d-2}{2d}}  \sm  \int_0^\infty \e^{-c_1s^{\ff 2 d} t^{-\bb} } s^{\ff{3d-2}{2d}}\d s\sm t^{\ff{\bb (5d-2)}4},\ \ t\ge 1.\end{align*}
 Similarly,  \eqref{CC'}  implies
\beg{align*} &t \|\nn f_{t,\bb}\|_\infty \sm \sum_{m=1}^\infty\ff{  \e^{-(\ll_m-\ll_0)t^{-\bb}} m^{\ff{3d+4}{2d} }} {(\ll_m-\ll_0)^2 }   \\
 &\sm   \sum_{m=1}^\infty   \e^{-c_1m^{\ff 2 d} t^{-\bb} } m^{\ff{3d-4}{2d} } \sm t^{\ff{\bb (5d-4)}4},\ \ t\ge 1.\end{align*}
Then the proof is finished.
 \end{proof}

\beg{lem}\label{LA5} For any $\bb\in (0,\ff 1 {20 d}]$, there exits a constant $c>0$ such that 
$$  t^2\W_2(\tt \mu_{t,\bb}^\nu,\mu_0)^2 \ge \ff {1-ct^{-1}} {\{\mu(\phi_0)\nu(\phi_0)\}^2} \sum_{m=1}^\infty \ff{\{\mu(h\phi_0)\mu_0(\phi_m)+ \mu(\phi_0) \nu(\phi_m)\}^2}{(\ll_m-\ll_0)^3}-ct^{-\ff 1 4}. $$
\end{lem}

\beg{proof}
To estimate $\W_2(\tt \mu_{t,\bb}^\nu,\mu_0)$ from below by using the argument in \cite{AMB,WZ20}, we take
$$\varphi_\theta^\vv:= -\vv \log P_{\ff{\vv \theta}2}^0 \e^{-\vv^{-1} f_{t,\bb}},\ \ \theta\in [0,1], \vv>0.$$
We have $\varphi_0^\vv=f_{t,\bb}$, $\|\varphi_\theta^\vv\|_\infty\le \|f_{t,\bb}\|_\infty$, and by \cite[Lemma 2.9]{WZ20},   there exists a constant $c_1>0$ such that  for any $\vv\in (0,1)$,
\beg{align*} &\varphi_1^\vv(y)- \varphi_0^\vv(x)\le \ff 1 2 \big\{\rr(x,y)^2 +\vv \|(L_0f_{t,\bb})^+\|_\infty + c_1\ss\vv \|\nn f_{t,\bb}\|_\infty^2\big\},\ \ x,y\in M,\\
&\int_M(\varphi_0^\vv-\varphi_1^\vv)\d\mu_0\le \ff 1 2 \int_M|\nn f_{t,\bb}|^2 \d\mu_0 +c_1 \vv^{-1} \|\nn f_{t,\bb}\|_\infty^4. \end{align*}
Therefore, by the Kantorovich dual formula, $\varphi_0^\vv=f_{t,\bb}$ and the integration by parts formula
$$\int_M  f_{t,\bb} \tt\rr_{t,\bb}^\nu\d\mu_0= \int_M f_{t,\bb} L_0 f_{t,\bb}\d\mu_0=-\int_M|\nn f_{t,\bb}|^2\d\mu_0,$$
we find a constant $c>0$ such that
\beq\label{OP} \beg{split} &c\big(\vv \|L_0 f_{t,\bb}\|_\infty +\vv^{\ff 1 2} \|\nn f_{t,\bb}\|_\infty^2\big) +\ff 1 2 \W_2(\tt\mu_{t,\bb}^\nu,\mu_0)^2 \ge \int_M \varphi_1^\vv\d\mu_0- \int_M \varphi_0^\vv \d\tt\mu_{t,\bb}^\nu\\
&=\int_M (\varphi_1^\vv-\varphi_0^\vv)\d\mu_0 -\int_M f_{t,\bb} \tt\rr_{t,\bb}^\nu \d\mu_0=\int_M (\varphi_1^\vv-\varphi_0^\vv)\d\mu_0 -\int_M f_{t,\bb} L_0 f_{t,\bb}  \d\mu_0 \\
&\ge  \ff  1 2 \int_M|\nn f_{t,\bb}|^2 \d\mu_0 -c \vv^{-1} \|\nn f_{t,\bb}\|_\infty^4.\end{split}\end{equation}
Taking $\vv= t^{-\ff 3 2}$ and applying Lemma \ref{L3.1}, when $\bb\le \ff 1{20d}$ we find a constant $c'>0$ such that
\beq\label{SPP} t^2\W_2(\tt\mu_{t,\bb}^\nu, \mu_0)^2 \ge t^2\mu_0(|\nn f_{t,\bb}|^2) - c't^{-\ff 1 4},\ \ t\ge t_0.\end{equation}
 Combining this with \eqref{GY3} and  \eqref{SPP},   we complete the proof.
\end{proof}

\beg{lem}\label{LA3}  There exist  constants $c,t_0>0$ such that for any $\nu=h\mu\in \scr P_0$ with $h\in L^2(\mu)$,  $\tt\mu_{t,\bb}^\nu, \tt\mu_t^\nu\in \scr P_0$ for $t\ge t_0$ and 
$$t\W_2(\tt\mu_{t,\bb}^\nu, \tt\mu_t^\nu)\le c \|h\|_{L^2(\mu)} t^{-\bb},\ \ t\ge t_0.$$
\end{lem}
\beg{proof}  $\tt\mu_{t,\bb}^\nu, \tt\mu_t^\nu\in \scr P_0$ for large $t$ is implied by Lemma \ref{LGY2}. Next, by \eqref{ALS}, we have
$$\scr M(\tt h_{t}^\nu,\tt h_{t,\bb}^\nu)\ge \ff 1 2,$$
so that  \cite[Lemma 2.3]{WZ20} implies
\beq\label{RM02} \W_2(\tt\mu_{t,\bb}^\nu, \tt\mu_{t}^\nu)^2\le \int_M\ff{|\nn L_0^{-1}(\tt h_{t}^\nu-\tt h_{t,\bb}^\nu)|^2}{\scr M(\tt h_{t}^\nu,\tt h_{t,\bb}^\nu)}\d\mu_0\le 2\mu_0(|\nn L_0^{-1}(\tt \rr_{t}^\nu-\tt \rr_{t,\bb}^\nu)|^2) .\end{equation}
To estimate the upper bound in this inequality, we first observe that by \eqref{GYY} and \eqref{GY2},  when $\nu=h\mu$ we have 
\beq\label{RM03'} \beg{split} &L_0^{-1} (\tt \rr_{t,\bb}^\nu-\tt \rr_{t}^\nu) =L_0^{-1} (P_{t^{-\bb}}^0 \tt \rr^\nu_{t}-\tt \rr^\nu_{t})=\int_0^{t^{-\bb}} P_r^0\tt\rr_t^\nu \d r\\
&=   \ff 1 {t\nu(\phi_0P_t^0\phi_0^{-1})} \int_0^{t^{-\bb}} (-L_0)^{-1}(P_r^0-\mu_0) g \, \d r, \end{split} \end{equation} where 
$$g:=  \mu(\phi_0) h\phi_0^{-1} +\nu(\phi_0)\phi_0^{-1}.$$ Since $\|h\|_{L^2(\mu)}\ge \mu(h)=1$, 
\beq\label{D*} \|g\|_{L^2(\mu_0)} \le\|\phi_0\|_\infty (1+ \|h\|_{L^2(\mu)})\le 2\|\phi_0\|_\infty \|h\|_{L^2(\mu)}.\end{equation} 
By      \eqref{000}, \eqref{D*}     and the fact
that  $(-L_0)^{-\ff 1 2}=c \int_0^\infty P_{s^2}^0\d s$ for some constant $c>0$,
we  find a constants $c_1,c_2>0$ such that 
 \beg{align*} &  \big\|\nn L_0^{-1}  (P_r^0-\mu_0)g\|_{L^2(\mu_0)}=\big\|L_0^{-\ff 1 2}  (P_r^0-\mu_0)g\|_{L^2(\mu_0)}\le \int_0^\infty \|(P_{r+s^2}^0-\mu_0)g\|_{L^2(\mu_0)} \d s\\
 &\le c_1  \|h\|_{L^2(\mu)} \int_1^\infty \e^{-(\ll_1-\ll_0)(s^2+r)} \d s\le c_2\|h\|_{L^2(\mu)},\ \ r\in [0,1].\end{align*}     Therefore, 
   by \eqref{GY3} and \eqref{RM03'}, we obtain
\beg{align*}  \|\nn L_0^{-1}(\tt\rr_{t,\bb}^\nu-\tt\rr_t^\nu) \|_{L^2(\mu_0)}  \sm \ff {1}t \int_0^{t^{-\bb} }\big\|\nn L_0^{-1}  (P_r^0-\mu_0)g\|_{L^2(\mu_0)}\d r
\sm t^{-(1+\bb)}\|h\|_{L^2(\mu)},\ \ t\ge t_0.\end{align*} 
   Combining this with \eqref{RM02} we finish the proof. 
 \end{proof}

We are now ready to prove the following result.

\beg{prp} For any $\nu\in \scr P_0$, 
\beq\label{JST}   \liminf_{t\to\infty} \big\{t^2\W_2(\mu_t^\nu,\mu_0)^2\big\}\ge I>0,\end{equation}
and  $I<\infty$ provided either $d\le 7$, or  $d\ge 7$ but $\nu=h\mu$ with  $  h\in L^{\ff{2d}{d+6}}.$    
 \end{prp}

\beg{proof}  Let $\bb\in (0,\ff 1 {20 d}]$. By \eqref{OK}, Lemma \ref{LA5} and Lemma \ref{LA3}, there exist    constants $c,t_0>0$ such that for $\nu= h\mu\in \scr P_0$ and $t\ge t_0$, 
\beg{align*}& t\W_2(\mu_{t}^\nu, \tt\mu_t^\nu)\le c\|h\|_{L^2(\mu)} t^{-\bb t},\\
&t \W_2(\tt\mu_{t,\bb}^\nu,\mu_0)\ge \big(\{(1-ct^{-1})I-ct^{-\ff 1 4})^+\big\}^{\ff 1 2},\\
&t\W_2(\mu_t^\nu, \tt\mu_t^\nu) \le ct \e^{-(\ll_1-\ll_0)t/2}\|h\|_{L^2(\mu)}^{\ff 1 2}.\end{align*}
Then 
\beq\label{MX} t\W_2(\mu_t^\nu,\mu_0)\ge  \big(\{(1-ct^{-1})I-ct^{-\ff 1 4})^+\big\}^{\ff 1 2}- c\|h\|_{L^2(\mu)} t^{-\bb t}-  ct \e^{-(\ll_1-\ll_0)t/2}\|h\|_{L^2(\mu)}^{\ff 1 2},\ \ t\ge t_0.\end{equation}
In general, let $\mu_{t,\vv}^\nu=\mu_{t-\vv}^{\nu_\vv}$ be as in the proof of Proposition \ref{PA}.  Applying  \eqref{MX} to $\mu_{t,t^{-2}}^\nu$ replacing $\mu_y^\nu$ and using \eqref{MX0}, \eqref{LMT}, 
we obtain
$$ \liminf_{t\to\infty} \big\{t\W_2(\mu_{t,t^{-2}}^\nu,\mu_0)\big\}\ge  \ss I,$$
which together with \eqref{*AG3} proves \eqref{JST}. 
 
  It remains to prove $I>0$ and $I<\infty $ the under given conditions, where due to \eqref{DMM}, $I<\infty$ is equivalent to
\beq\label{II0}\beg{split} &I':= \sum_{m=1}^{\infty}\ff{\nu(\phi_m )^2}{\ll_m^3}<\infty.\end{split}\end{equation}
Below we first prove $I>0$ then show $I'<\infty$ under the given conditions.

(a) $I>0$.  If this is not true, then 
 $$ \mu(h\phi_0)\mu(\phi_m) = - \mu(\phi_0) \mu(h\phi_m),\ \ m\ge 1.$$ Combining this with the representation in $L^2(\mu)$
 $$f= \sum_{m=0}^\infty \mu(f\phi_m)\phi_m,\ \ f\in L^2(\mu),$$
 where the equation holds point-wisely if $f\in C_b(M)$ by the continuity, we obtain
 \beg{align*} &\mu(\phi_0) \nu(f)= \sum_{m=0}^\infty \mu(f\phi_m) \mu(\phi_0)\nu(\phi_m) = 2 \mu(f\phi_0)\nu(\phi_0)\mu(\phi_0) - \sum_{m=0}^\infty \mu(f\phi_m) \mu(\phi_m)\nu(\phi_0)\\
 &=  2 \mu(f\phi_0)\nu(\phi_0)\mu(\phi_0) -  \nu(\phi_0)\mu(f),\ \ \ f\in C_b(M).\end{align*}
 Consequently,
 $$0\le \mu(\phi_0)\ff{\d \nu}{\d\mu}= 2\phi_0\nu(\phi_0)\mu(\phi_0) -  \nu(\phi_0),  $$ which is however impossible since the upper bound is negative in a neighborhood of $\pp M$, because   $\nu(M^\circ)>0$ implies $\nu(\phi_0)>0$ for $\phi_0>0$ in $M^\circ$,
 and $\phi_0$ is continuous with $\phi_0|_{\pp M}=0.$  Therefore, we must have $I>0.$

(b)   $I'<\infty.$    Let $\{h_n\}_{n\ge 1}$ be a sequence of probability density functions with respect to $\mu$ such that 
\beq\label{NUM} \nu_n:=h_n\mu\to \nu \ \text{ weakly\ as \ }n\to\infty.\end{equation}
By   the spectral representation for $(-L)^{-\ff 3 2}$,  and applying the Sobolev inequality \eqref{SB} with $p= \ff{2d}{d+6}\lor 1,$ we obtain 
 \beq\label{IN} I_n':= \sum_{m=1}^\infty \ff{\nu_n(\phi_m)^2} {\ll_m^3} \le \|(-L)^{-\ff 3 2} h_n\|_{L^2(\mu)}^2 \le K^2 \|h_n\|_{L^{\ff{2d}{d+6}\lor 1}(\mu)}^2,\ \ n\ge 1.\end{equation}
It is easy to see that for $d\le 6$ we have $\ff{2d}{d+6}\le 1$, so that $ \|h_n\|_{L^{\ff{2d}{d+6}\lor 1}(\mu)}=\mu(h_n)=1.$ Combining this with \eqref{NUM}, \eqref{IN}  and applying Fatou's lemma, we derive 
$$I' = \sum_{m=1}^\infty \liminf_{n\to\infty} \ff{\nu_n(\phi_m)^2} {\ll_m^3} \le \liminf_{n\to\infty} I_n'\le K^2<\infty,\ \ d\le 6.$$
 Finally, when     $d\ge 7$ and $\nu=h\mu$ with  $  h\in L^{\ff{2d}{d+6}}(\mu),$  by applying \eqref{IN} to $h_n=h$ we prove $I'<\infty$. 
    \end{proof}

\paragraph{Acknowledgement.} The author  would like to thank the referee for helpful comments and corrections.

\end{document}